\documentclass[letterpaper, 10 pt, conference]{ieeeconf}  

\IEEEoverridecommandlockouts                              

\usepackage{amsmath}
\usepackage{amssymb}
\usepackage{algorithm}
\usepackage{algpseudocode}

\usepackage[caption=false,font=footnotesize]{subfig}

\usepackage{multirow}
\usepackage[noadjust]{cite}

\usepackage{graphicx} 

\title{\LARGE \bf
Designing Real-Time Prices to Reduce Load Variability with HVAC}

\author{John Audie Cabrera$^{2}$, Yonatan Mintz$^{1}$, Jhoanna Rhodette Pedrasa$^{2}$, and Anil Aswani$^{1}$
\thanks{*This work was supported in part by the Philippine-California Advanced Research Institutes (PCARI) and NSF Award CMMI-1450963.}
\thanks{$^{1}$Y. Mintz and A. Aswani are with the Department of Industrial Engineering and Operations Research, University of California, Berkeley, CA 94720 USA 
        {\tt\small ymintz@berkeley.edu, aaswani@berkeley.edu}}%
\thanks{$^{2}$J.A. Cabrera and J.R. Pedrasa are with the Electrical and Electronics Engineering Institute, University of the Philippines, Diliman, Quezon City, Philippines 1101
        {\tt\small john\_audie.cabrera@upd.edu.ph, jipedrasa@up.edu.ph}}%
}

\begin{document}

\maketitle
\thispagestyle{empty}
\pagestyle{empty}

\begin{abstract} 
Utilities use demand response to shift or reduce electricity usage of flexible loads, to better match electricity demand to power generation.  A common mechanism is peak pricing (PP), where consumers pay reduced (increased) prices for electricity during periods of low (high) demand, and its simplicity allows consumers to understand how their consumption affects costs.  However, new consumer technologies like internet-connected smart thermostats simplify real-time pricing (RP), because such devices can automate the tradeoff between costs and consumption.  These devices enable consumer choice under RP by abstracting this tradeoff into a question of quality of service (e.g., comfort) versus price.  This paper uses a principal-agent framework to design PP and RP rates for heating, ventilation, and air-conditioning (HVAC) to address adverse selection due to variations in consumer comfort preferences.  We formulate the pricing problem as a stochastic bilevel program, and numerically solve it by reformulation as a mixed integer program (MIP).  Last, we compare the effectiveness of different pricing schemes on reductions of peak load or load variability.  We find that PP pricing induces HVAC consumption to spike high (before), spike low (during), and spike high (after) the PP event, whereas RP achieves reductions in peak loads and load variability while preventing large spikes in electricity usage.
\end{abstract}

\section{Introduction}
High demand variability stresses the electrical grid by increasing the mismatch with supply, and it is costly for utilities because it requires adding redundant power generation.  Demand response is an alternative that induces consumers to reduce or shift their consumption by setting prices by time of day \cite{Newsham2011,Gyamfi2012,Gyamfi2013,Sun2013,Strbac2008}.  For example, peak pricing (PP) reduces the peak demand of electricity by charging consumers reduced (increased) rates for electricity during periods of low (high) demand.  This is a common structure for demand response programs because the simplicity of PP allows consumers to understand how their consumption impacts their costs.

Real-time pricing (RP) of electricity is less common because historically the complex pricing structure of RP makes it difficult for consumers to match consumption to prices.  However, new consumer technologies like internet-connected smart thermostats \cite{AswaniMasterTanejaEtAl2012,AswaniMasterTanejaEtAl2012b,MaasoumyRosenbergSangiovanni-VincentelliEtAl2014,ZugnoMoralesPinsonEtAl2013,BorscheOldewurtelAndersson2013,vrettos2013predictive} simplify RP, because such devices can automate the tradeoff between costs and consumption.  These devices simplify RP by abstracting this tradeoff into a question of quality of service (e.g., comfort) versus price, which is easier for consumers to understand.

This paper designs PP and RP electricity rates using realistic, validated models of heating, ventilation, and air-conditioning (HVAC) \cite{AswaniMasterTanejaEtAl2012,AswaniMasterTanejaEtAl2012a}, and there are three contributions.  First, we use a principal-agent model \cite{LaffontMartimort2009,aswani2012incentive} to formulate the problem of a utility designing rates for HVAC that responds to prices, where the consumer has an acceptable (but unknown to the utility) comfort level.  The challenge is that prices must be designed so that inflexible (with respect to comfort) consumers do not get excessive benefits relative to flexible consumers, since flexible consumers provide more benefits to the utility.  Second, we pose the design problem as a mixed integer program (MIP).  Third, we present numerically solvable approximations of this MIP, and then evaluate the impact of the resulting PP and RP rates.

\subsection{PP for HVAC Demand Response}

HVAC is arguably the most significant target for demand response since it the largest source of energy consumption in most buildings \cite{afram2014}.  This is relevant from the standpoint of utilities because HVAC use is obviously correlated with high outdoor temperature, which means that HVAC usage in different buildings is strongly correlated with each other and is an important contributor to peak demand \cite{mendoza2012}.  As a result, many studies have considered different aspects of PP for demand response of HVAC.  A large number of demand response programs that have been implemented by utility companies use PP to reduce peak load \cite{Newsham2011,Gyamfi2012,Gyamfi2013,Sun2013,Strbac2008}, and such programs have been found to provide varying levels of value to utilities.  Within the controls literature, the use of model predictive control (MPC) techniques is particularly popular for demand response of HVAC \cite{kelman2011,parisio2014,mintz2016,mintz2017behavioral} because of the ability of MPC to handle complex constraints.

\subsection{RP for HVAC Demand Response}

Recent work studied RP design for HVAC that automates price-responsiveness.  One approach uses stochastic differential equations to design prices \cite{YangCallawayTomlin2014,YangCallawayTomlin2015}, and this work found a benefit to RP for a simplified HVAC model.  In contrast, we consider in this paper the rate design problem using realistic, validated models of HVAC \cite{AswaniMasterTanejaEtAl2012,AswaniMasterTanejaEtAl2012a}.  Another body of work \cite{avci2012residential,avci2013model} considers RP design using realistic HVAC models.  Our paper differs in two substantive ways.  The first is we use a different notion of comfort: Comfort in \cite{avci2012residential,avci2013model} was defined using the temperature set-point, whereas in our paper we define comfort using allowable deviations in the temperature from the desired value.  The second is we consider adverse selection, which are the issues caused when an inflexible (with respect to comfort) consumer accepts a rate designed for a flexible consumer, in our rate design.

\subsection{Outline}

Sect. \ref{sec:cm} describes our model for the consumer and our model for the electric utility company, including the principal-agent model the utility uses to design the electricity rates.  The key feature of the model is the fact that consumers are either flexible or inflexible with regards to their comfort, but this information is hidden from the electric utility.  The electricity rate will not be efficient for the utility if it does not account for this information asymmetry (formally known as \emph{adverse selection}).  Next, Sect. \ref{sec:nspp} describes how to numerically solve the rate design problem using an MIP reformulation of the principal-agent model.  As part of our approach, we derive relaxations that facilitate fast numerical solution.  We conclude with Sect. \ref{sec:nr}, which numerically solves the pricing problem and then compares the impact of PP and RP on electricity consumption by HVAC.

\section{Model of Consumer and Electric Utility}
\label{sec:cm}

In this section, we present our model for the consumer and the electric utility.  We also formally define the problem of using a principal-agent framework to design either PP or RP electricity prices for HVAC demand response.

\subsection{Consumer Model}

The first part of our model defines comfort in relation to deviations in room temperature from the desired value: Consumers are inflexible ($\pm 2^\circ$C deviation from desired temperature) or flexible ($\pm 3^\circ$C deviation from desired temperature) in their comfort, and these ranges are from the ASHRAE 55 standard \cite{Standard2010} that defines quantitative models of occupant comfort.  We use $T_d$ to refer to a consumer's desired room temperature, and the $\overline{T},\underline{T}$ are the upper and lower bounds of comfort for the consumer.  So if the consumer is inflexible, then $\underline{T} = T_d-2$ and $\overline{T} = T_d+2$ .  Similarly, if the consumer is flexible, then $\underline{T} = T_d-3$ and $\overline{T} = T_d+3$.

The next part of our model describe the room temperature dynamics and provides an energy model for the consumer.  We use a linear time-invariant model for room temperature
\begin{equation}
T_{n+1} =k_{r}T_{n}+k_{c}u_{n}+k_{w}w_{n}+q_n,\label{eq: dynamic}
\end{equation}
where $T_{n}$, $u_{n}$, $w_{n}$, $q_n$ are room temperature, HVAC control input, outside temperature, and heating load due to occupancy, respectively, and each time step is a 15 min interval.  This model has been validated \cite{AswaniMasterTanejaEtAl2012,AswaniMasterTanejaEtAl2012a}.  The total energy usage of the consumer is $\sum_{n=1}^N(b_n + pu_n)$, where $b_n$ is nondeferrable electricity load, $p$ is a constant that converts input $u_n$ to energy consumption \cite{AswaniMasterTanejaEtAl2012,AswaniMasterTanejaEtAl2012a}, and $N$ is a horizon.

An important component of our model characterizes the HVAC controller, which automates the tradeoff between room temperature and electricity consumption.  In particular, we assume that the HVAC is controlled by MPC:
\begin{equation}
\label{eqn:mpc}
\begin{aligned}
\min\ & \textstyle\sum_{n=1}^{N}\big((T_{n}-T_{d})^{2}+\gamma c_{n}u_{n}\big)\\
\mathrm{s.t.}\ & T_{n+1}=k_{r}T_{n}+k_{c}u_{n}+k_{w}w_{n}+q_n\\
& T_n \in [\underline{T},\overline{T}], u_n \in [0,\overline{u}],\quad \text{ for } n = 1,\ldots,N
\end{aligned}
\end{equation}
where $\gamma$ is a constant that trades off temperature and electricity usage, $\overline{u}$ is the maximum control input, and $c_n$ is the price of electricity at time $n$.

The last part of the model describes what information is known by the consumer (and implicitly known by the HVAC controller).  The variable
\begin{multline}
\theta = \big\{k_r,k_c,k_w,w_n,q_n,b_n,\gamma,T_d,\underline{T},\overline{T},\overline{u}, \\\text{for } n=1,\ldots,N\big\}
\end{multline}
completely characterizes each consumer, and it is known as \emph{type} in the principal-agent literature \cite{LaffontMartimort2009,aswani2012incentive}.  (The value $p$ is a constant known by everyone.)  We assume that the consumer (and HVAC controller) exactly knows the value of $\theta$, and knows the electricity price $\textbf{c} = \{c_1,\ldots,c_N\}$.  Moreover, we use $J(\textbf{c}; \theta)$ to refer to the minimum value of (\ref{eqn:mpc}), and $u^*(\textbf{c}; \theta)$ refers to the minimizer of (\ref{eqn:mpc}).

\subsection{Model of Electric Utility Company}

An important component in the electric utility model is the information asymmetry between the utility and consumers.  Specifically, we assume the utility does not know $\theta$ for any single customer.  Instead, the utility knows the overall probability distribution for $\theta$.  (Recall the utility and consumers know $p$, which is a constant.)  We also assume that both the utility and consumers know the electricity price $\textbf{c}$.

The next element in the utility model describes the goal of the electricity pricing for demand response.  If the goal is to reduce peak load, then the utility aims to minimize
\begin{equation}
V_{p} = \textstyle\mathbb{E}_\theta\Big(\sum_{n = t_1}^{t_2} u^*_n(\textbf{c}; \theta)\Big),
\end{equation}
where $[t_1,t_2]$ is a time range during which the peak load is anticipated by the utility. If the goals is to reduce load variability, then the utility aims to minimize
\begin{equation}
V_{l} = \textstyle\mathbb{E}_\theta\Big(\mathrm{var}_n\big(b_n + u_n^*(\textbf{c}; \theta)\big)\Big),
\end{equation}
where $\mathrm{var}_n(\cdot)$ is the variance over $n=1,\ldots,N$.  We will consider designing PP and RP for both goals.

The electric utility is interested in designing $\textbf{c}$, and we describe the constraints that characterize PP and RP rates.  If the utility is designing PP rates, then this means they are selecting from
\begin{equation}
\mathcal{C}_{pp} = \left\{\textbf{c} : 
\begin{aligned}
&c_n = c_{t_1}, &\text{ for } n \in[t_1,t_2]\\
& c_n = c_1, &\text{ for } n \in \{1,\ldots,N\}\setminus[t_1,t_2]
\end{aligned}\right\}.
\end{equation}
This expresses prices that are constant within the peak period $[t_1,t_2]$, and constant (with a possibly different value) outside of the peak period.  Similarly, if the utility is designing RP rates, then this means they are selecting from
\begin{equation}
\mathcal{C}_{rp} = \left\{\textbf{c} : 
\begin{aligned}
&c_1 = c_N\\
&|c_{n+1}-c_n| \leq \rho, &\text{ for } n =1,\ldots,N_1
\end{aligned}\right\}.
\end{equation}
This expresses prices that are equal at the beginning and end of the horizon, and such that the rate of change is bounded by a constant $\rho$.  Lastly, we use $\textbf{f} = \{f,\ldots,f\}$ to refer to a flat pricing structure, and $f$ in particular refers to the existing electricity price prior to the introduction of the demand response pricing.

\subsection{Principle-Agent Model for Pricing}

\begin{table}
\centering
\begin{tabular}{ccccc}
 & $k_{r}$ & $k_{c}$ & $k_{w}$ & average $q_n$\tabularnewline
\hline 
\hline 
Room 1 & 0.63 & 2.64 & 0.10 & 6.78\\
Room 2 & 0.43 & 1.95 & 0.18 & 9.44\\
\end{tabular}
\caption{\label{tabel:thermal_coeff}Temperature Model Coefficients}
\end{table}

The last part of the model for the utility describes the principal-agent formulation used to design electricity prices.  In particular, we assume the utility solves
\begin{equation}
\label{eqn:pam}
\begin{aligned}
\min\ & \textstyle V + \lambda\cdot\mathbb{E}_\theta\Big(\sum_{n=1}^N\big(f_n^{\vphantom{*}}u_n^*(\mathbf{f}; \theta) - c_n^{\vphantom{*}}u_n^*(\mathbf{c}; \theta)\big)\Big) \\
\mathrm{s.t}\ & J(\textbf{c}; \theta) \leq J(\textbf{f}; \theta)\\
&\mathbf{c}\in\mathcal{C}\\
& c_n \in [\underline{c},\overline{c}],\quad \text{for } n=1,\ldots,N 
\end{aligned}
\end{equation}
to design the electricity rates, where $V$ is either $V_{p}$ (to minimize peak load) or $V_{v}$ (to minimize load variance), and $\mathcal{C}$ is either $\mathcal{C}_{pp}$ (for PP) or $\mathcal{C}_{rp}$ (for RP).  Note the $\underline{c},\overline{c}$ are bounds on the minimum and maximum electricity rate, respectively.

Here, $\sum_{n=1}^N\big(f_nu_n^*(\mathbf{f}; \theta) - c_nu_n^*(\mathbf{c}; \theta)\big)$ is the amount of revenue the utility loses from implementing the new pricing $\textbf{c}$ (relative to the existing rate $\mathbf{f}$), and so this means $\lambda$ is a constant that the utility uses to tradeoff achieving the demand response goal with revenue loss.  We do not include the nondeferrable electricity load $b_n$ when defining revenue loss, because in our setting the electricity rates for the nondeferrable electricity load are different (and left unchanged) from the rates $\textbf{c}$ for HVAC electricity consumption.  

There are two game-theoretic considerations that must be discussed when defining and solving principal-agent models \cite{LaffontMartimort2009,aswani2012incentive}. The constraint $J(\textbf{c}; \theta) \leq J(\textbf{f}; \theta)$ is known as a \emph{participation constraint}, and it ensures that the new electricity rates $\textbf{c}$ are such that the overall utility of the consumer under the new rates $\mathbf{c}$ is equal or better than the overall utility of the consumer under the original rate $\mathbf{f}$.  The second game-theoretic aspect to be discussed is adverse selection.  We mitigate adverse selection by minimizing the expectation (with respect to type $\theta$) of the goal $V$ and revenue loss.

\section{Numerical Solution of Pricing Problem}

\label{sec:nspp}

This section studies how to solve the principal-agent model (\ref{eqn:pam}).  The main difficulty is that (\ref{eqn:pam}) is a bilevel program \cite{ColsonMarcotteSavard2007,aswani2016duality}, which means that (\ref{eqn:pam}) is an optimization problem in which some variables are solutions to optimization problems themselves.  In particular, recall that $u^*(\textbf{c}; \theta)$ is the minimizer to (\ref{eqn:mpc}).  In order to solve (\ref{eqn:pam}), we first show how the problem can be reformulated as a MIP.  Then we describe some relaxations that facilitate numerical solution of the MIP.

\subsection{MIP Reformulation of Pricing Problem}

They key idea in reformulating (\ref{eqn:pam}) is to replace the convex optimization problem (\ref{eqn:mpc}) by the KKT conditions, which provides constraints that $u^*(\textbf{c}; \theta)$ must satisfy.  More specifically, the KKT conditions for (\ref{eqn:mpc}) can be written as the following set of mixed integer linear constraints:
\begin{equation}
\begin{aligned}
&T_{n+1}=k_{r}^{\vphantom{*}}T_{n}^{\vphantom{*}}+k_{c}^{\vphantom{*}}u_{n}^*(\textbf{c}; \theta) +k_{w}^{\vphantom{*}}w_{n}^{\vphantom{*}}+q_n^{\vphantom{*}}\\
& \gamma c_{n}-k_{c}\nu_{n}+\overline{\mu}_{n}-\underline{\mu}_{n}=0\\
&0\leq\overline{\mu}_{n}\leq M\eta_n,\\
&0\leq\underline{\mu}_{n}\leq M\zeta_n\\
&\overline{u}\eta_n+\underline{u}\left(1-\eta_n\right)\leq u_{n}^*(\textbf{c}; \theta)\leq\underline{u}\zeta_n+\overline{u}\left(1-\zeta_n\right)\\
&\eta_n, \zeta_n \in\{0,1\}, \quad \text{for } i = 1,\ldots,N-1
\end{aligned}
\end{equation}
and also that
\begin{equation}
\begin{aligned}
&(T_{n}-T_{d})+\nu_{n-1}-k_{r}\nu_{n}+\overline{\xi}_{n}-\underline{\xi}_{n}=0\\
&0\leq\overline{\xi}_{n}\leq Mx_{n}\\
&0\leq\underline{\xi}_{n}\leq My_{n}\\
&\overline{T}x_{n}+\underline{T}\left(1-x_{n}\right)\leq T_{n}\leq\underline{T}y_{n}+\overline{T}\left(1-y_{n}\right)\\
&x_n, y_n \in\{0,1\}, \quad \text{for } 2 = 1,\ldots,N
\end{aligned}
\end{equation}
where $M > 0$ is a sufficiently large constant \cite{Fortuny-AmatMcCarl1981}. 

The problem (\ref{eqn:pam}) becomes an infinite dimensional MIP, after a few more reformulations.  The first is to observe that $\mathbb{E}_\theta(f_n^{\vphantom{*}}u_n^*(\mathbf{f}; \theta))$ is a constant, and so can be removed from the objective function.  The second is to note that $J(\textbf{f}; \theta)$ is also a constant since it does not depend on any decision variables.  The third reformulation is to substitute $J(\textbf{c}; \theta)$ with $\sum_{n=1}^{N}\big((T_{n}^{\vphantom{*}}-T_{d^{\vphantom{*}}})^{2}+\gamma c_{n}^{\vphantom{*}}u_{n}^*(\textbf{c}; \theta)\big)$.  Though this yields an infinite dimensional problem, using sample average approximation (SAA) \cite{kleywegt2002sample,wang2008sample} to approximate the reformulation gives a finite dimensional MIP.

\subsection{Relaxation of Pricing Problem}

The reformulated MIP described above is still difficult to solve because it involves nonconvex quadratic terms $c_{n}^{\vphantom{*}}u_{n}^*(\textbf{c}; \theta)$, and so additional relaxations are needed so that the price design problem can be solved using standard numerical optimization software.  The quadratic term is relaxed using the McCormick envelope \cite{McCormick1976} to
\begin{equation}
\begin{aligned}
r_n \geq \underline{c}u_{n}^*(\textbf{c}; \theta) + \underline{u}c_n - \underline{u}\cdot\underline{c}\\
r_n \geq \overline{c}u_{n}^*(\textbf{c}; \theta) + \overline{u}c_n - \overline{c}\cdot\overline{u}\\
r_n \leq \overline{c}u_{n}^*(\textbf{c}; \theta) + \underline{u}c_n - \overline{c}\cdot\underline{u}\\
r_n \leq \underline{c}u_{n}^*(\textbf{c}; \theta) + \overline{u}c_n - \underline{c}\cdot\overline{u}
\end{aligned}
\end{equation}
for $n = 1,\ldots,N$.  With this relaxation, the SAA form of the reformulated problem is a mixed-integer quadratic program (MIQP), which can be solved using existing software.

However, numerical solution of MIQP's can be slow.  So we next describe two additional relaxations that speed up computation by approximating the MIQP using a mixed-integer linear program (MILP), which can typically be numerically solved faster.  First, we replace $(T_n-T_d)^2$ with $4|T_n-T_d|^2$ since $(T_n-T_d) \leq 3|T_n-T_d|^2$ when $|T_n-T_d| \leq 3$ as is the case from our assumptions about comfort.  Second, we replace $\mathrm{var}_n\big(b_n + u_n^*(\textbf{c}; \theta)\big)$ with $N^{-1}\sum_{n=1}^N|b_n + u_n^*(\textbf{c}; \theta) - m(\theta)|$, where $m(\theta) = \frac{1}{N}\sum_{n=1}^N u_n^*(\textbf{f}; \theta)$.  The idea is we approximate the variance by (a) replacing squares with absolute value, and (b) replacing the mean in the variance $\frac{1}{N}\sum_{n=1}^N u_n^*(\textbf{c}; \theta)$ with the mean $m(\theta)$.

\begin{table}
\centering
\begin{tabular}{llccc}
&&\textbf{Flat Rate} & \textbf{PP Rate} & \textbf{RP Rate}\\
\hline
\hline
\multirow{2}{*}{\textbf{Inflexible}}&Peak Load&28.3&27.0&27.6\\
&Load Variance & 0.49&0.54&0.42\\
\hline
\multirow{2}{*}{\textbf{Flexible}}&Peak Load&19.1&15.3&17.5\\
&Load Variance & 0.25 & 0.28 & 0.17\\
\end{tabular}
\caption{Pricing to Reduce Peak Load\label{tab:rpl}}
\end{table}

\section{Numerical Results}

\label{sec:nr}

In this section, we numerically solve our MILP relaxation of the pricing problem for a 24 hour horizon.  All of the calculations where conducted on laptop computer with dual core 2.5GHz processor and 8GB RAM using MATLAB with the CVX toolbox \cite{cvx} and the Gurobi solver \cite{GurobiOptimization2016}.  We finish by evaluating the quality of the designed electricity rates, and the results are summarized in Tables \ref{tab:rpl} and \ref{tab:rlv}.

\subsection{Values of Type Parameters}

For scenarios with PP and peak load reduction, we set the peak times to be 1pm--4pm. Our bounds on the electricity cost were $7\text{PhP}\leq c_{n}\leq20\text{PhP}$, where PhP is Philippines Pesos.  Parameters in the room temperature dynamics (\ref{eq: dynamic}) were chosen by uniformly sampling from the paramters in Table \ref{tabel:thermal_coeff}.  The first set of parameters are from \cite{AswaniMasterTanejaEtAl2012,AswaniMasterTanejaEtAl2012a}, while the second set of parameters were replicated using the same methodology from \cite{AswaniMasterTanejaEtAl2012,AswaniMasterTanejaEtAl2012a} with data from our UP-BRITE testbed located at the University of the Phillipines, Diliman.  We set the probability of a consumer to have high flexibility to be 0.2.  Scenario generation for outside temperature was performed using data from Weather Underground \cite{wunderground}, scenario generation for heating load due to occupancy was based on occupancy models, and scenario generation for nondeferrable electricity load was based on data from \cite{nrel}.

\subsection{Results and Discussion for PP}

\begin{figure*}
\centering
\includegraphics[trim={0.7in 0in 0.5in 0in},clip=true]{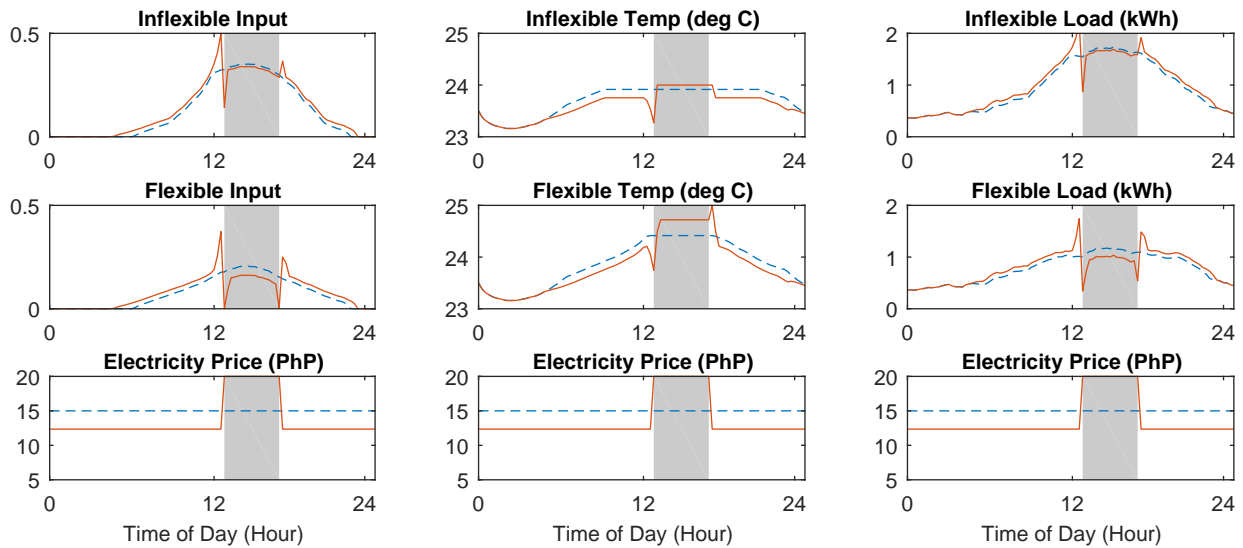}
\caption{Comparison of expected load response under constant pricing (dashed line) and PP designed to reduce peak load (solid line).\label{fig:pppl}} 
\end{figure*}

\begin{figure*}
\centering
\includegraphics[trim={0.7in 0in 0.5in 0in},clip=true]{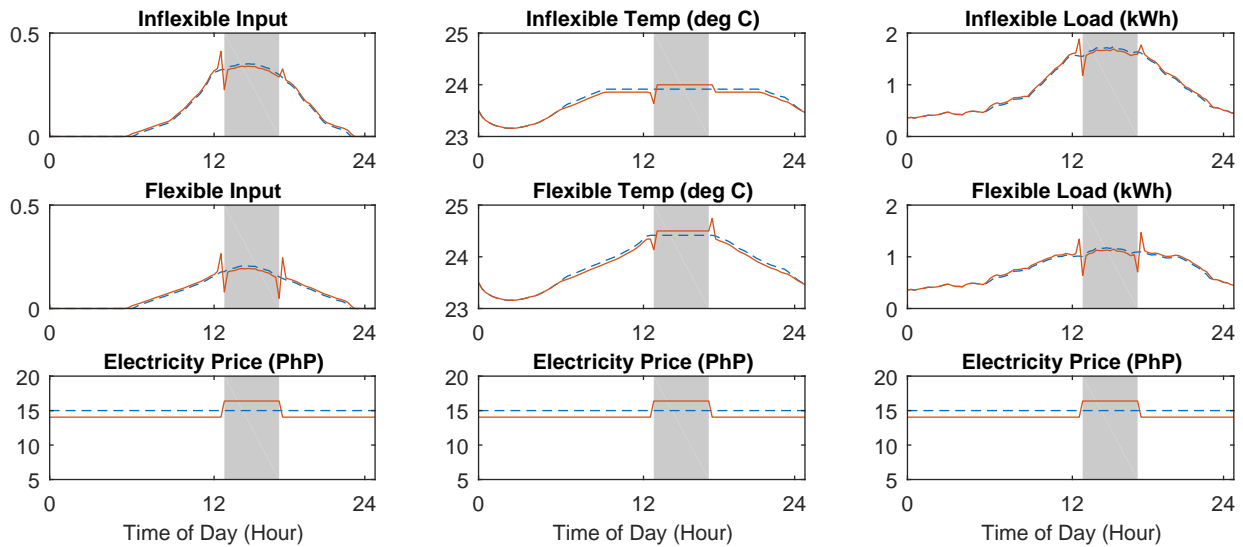}
\caption{Comparison of expected load response under constant pricing (dashed line) and PP designed to reduce load variance (solid line).\label{fig:ppv}} 
\end{figure*}

Results for PP for peak load reduction are shown in Fig. \ref{fig:pppl}.  PP is effective in reducing the peak load for both the flexible and inflexible consumers; but there is a side effect in which the HVAC has sharp increases in electricity consumption both prior to and after the peak period, as well as a sharp decrease in consumption at the start and end of the peak period.  This substantially increases the variability of the load profile.  Results for PP for load variance reduction are shown in Fig. \ref{fig:ppv}.  PP is not effective in decreasing load variability because sharp changes in electricity price induce the HVAC to make sharp changes in consumption.

\subsection{Results and Discussion for RP}

\begin{figure*}
\centering
\includegraphics[trim={0.7in 0in 0.5in 0in},clip=true]{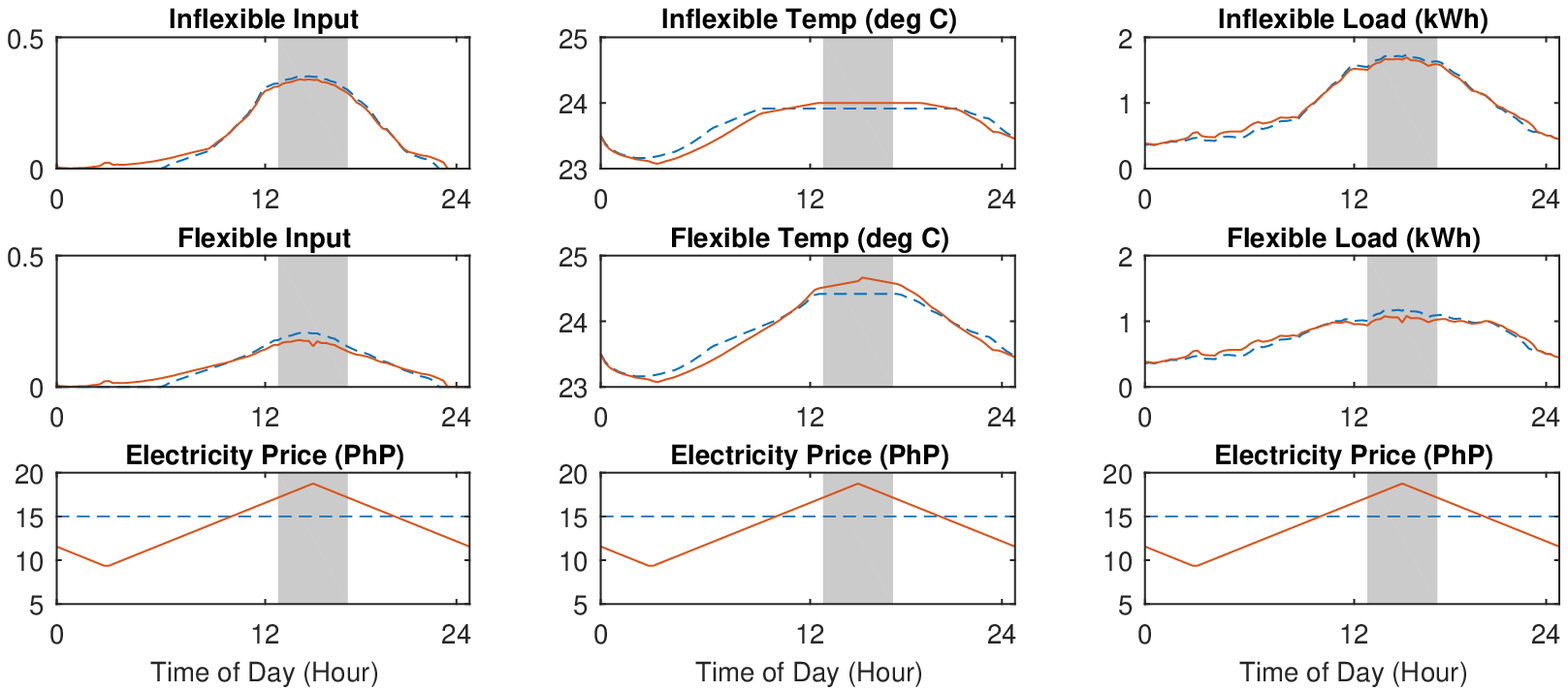}
\caption{Comparison of expected load response under constant pricing (dashed line) and RP designed to reduce peak load (solid line).\label{fig:rppl}} 
\end{figure*}

\begin{figure*}
\centering
\includegraphics[trim={0.7in 0in 0.5in 0in},clip=true]{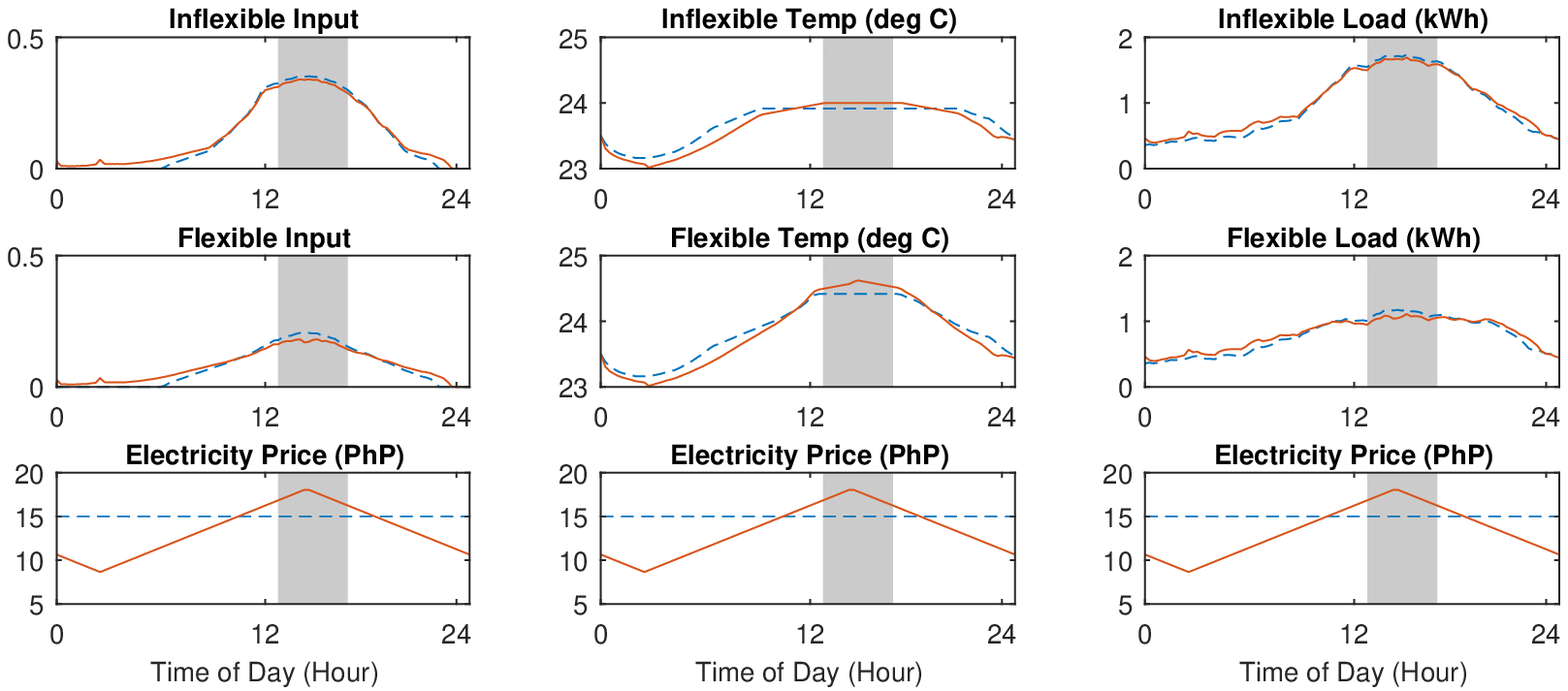}
\caption{Comparison of expected load response under constant pricing (dashed line) and RP designed to reduce load variance (solid line).\label{fig:rpv}} 
\end{figure*}

The results for RP for peak load reduction are shown in Fig. \ref{fig:rppl}.  The RP is effective in reducing the peak load for both the flexible and inflexible consumers, and it in fact also reduces the variance of the electricity load.  The results for RP for load variance reduction are shown in Fig. \ref{fig:rpv}.  The RP is effective in decreasing the variability of the total electricity load, and it also reduces the peak load for both the flexible and inflexible consumers.  The variance in load under this latter contract is lower than the variance under the former contract, but the difference is small.

\section{Conclusion}

\begin{table}
\centering
\begin{tabular}{llccc}
&&\textbf{Flat Rate} & \textbf{PP Rate} & \textbf{RP Rate}\\
\hline
\hline
\multirow{2}{*}{\textbf{Inflexible}}&Peak Load&28.3&27.3&27.6\\
&Load Variance & 0.49&0.49&0.41\\
\hline
\multirow{2}{*}{\textbf{Flexible}}&Peak Load&19.1&17.7&17.8\\
&Load Variance & 0.25 & 0.26 & 0.17\\
\end{tabular}
\caption{Pricing to Reduce Load Variance\label{tab:rlv}}
\end{table}

We studied the problem of designing PP and RP electricity rates using realistic, validated models of HVAC.  We used a principal-agent model to formulate the problem of a utility designing rates for HVAC that responds to prices, where the consumer has an acceptable (but unknown to the utility) comfort level.  We showed how this problem could be posed as numerically tractable MILP's, and then solved these MILP's to compare the efficacy of different pricing schemes.  We found that RP was substantially better at reducing load variability than PP, whereas PP was superior in reducing peak load.  Directions for future work include incorporating more detailed consumer models to better understand best practices for the design of incentives for effective demand response.

\bibliographystyle{IEEEtran}
\bibliography{IEEEabrv,hvar}

\end{document}